\documentclass[12pt, a4paper]{article}
\usepackage{}

\usepackage{amsfonts}
\usepackage{amsfonts}
\usepackage{bbm}
\usepackage{mathrsfs}
\usepackage{latexsym}
\usepackage{graphicx}
\usepackage{amssymb}

\usepackage{fullpage}
\usepackage{amssymb}
\usepackage{amsmath}
\usepackage{amsthm}

\newtheorem{theorem}{Theorem}[section]
\newtheorem{lemma}[theorem]{Lemma}

\newtheorem{conjecture}[theorem]{Conjecture}
\newtheorem{remark}[theorem]{Remark}

\title{{\bf A lower bound of the least signless Laplacian eigenvalue of a graph\thanks{ Supported by NSFC
(Nos. 11171290,~ 11271315, ~11201417, 11371195.)}~}}

\author{ Shu-Guang Guo$^a$\thanks{ E-mail addresses:
ychgsg@163.com (Guo), ygchen@njnu.edu.cn (Chen), yglong01@163.com
(Yu).}, ~  Yong-Gao Chen$^b$,  Guanglong Yu$^a$
 \\
{\footnotesize $^a$Department of Mathematics, Yancheng Teachers
University,}\\ {\footnotesize  Yancheng, 224002, Jiangsu, P.R.
China}\\
\footnotesize  $^b$School of Mathematical Sciences and Institute
of Mathematics,
\\  \footnotesize  Nanjing Normal University,  Nanjing
210023, P. R. CHINA }
\date{}

\begin{document}
\maketitle

\begin{abstract}
Let $G$ be a simple connected graph on $n$ vertices and $m$ edges.
In [Linear Algebra Appl. 435 (2011) 2570-2584], Lima et al. posed
the following conjecture on the least eigenvalue $q_n(G)$ of the
signless Laplacian of $G$: $\displaystyle  q_n(G)\ge
{2m}/{(n-1)}-n+2$. In this paper we prove a stronger result: For
any graph with $n$ vertices and $m$ edges, we have $\displaystyle
q_n(G)\ge {2m}/{(n-2)}-n+1 (n\ge 6)$.

\bigskip
\noindent {\bf AMS Classification:} 05C50

\noindent {\bf Keywords:} Graph;  Signless Laplacian; Least
eigenvalue
\end{abstract}
\baselineskip 20pt

\maketitle

\baselineskip 20pt

\section{Introduction}

For a undirected and simple graph $G$, let
 $A(G)$ be its adjacency matrix and let
$D(G)$ be the diagonal matrix of its degrees. The matrix
$Q(G)=D(G)+ A(G)$ is called the signless Laplacian matrix or
$Q$-matrix of $G$.  As usual, we shall index the eigenvalues of
$Q(G)$ in nonincreasing order and denote them as $q_1(G)\ge q_2(G)
\ge \cdots \ge q_n(G)\ge  0$.

The signless Laplacian eigenvalues of a graph have recently
attracted more and more researchers' attention, see
\cite{Cvetkovic} and its references. Lima et al. \cite{LOA}
surveyed some known results, presented some new ones about the
least signless Laplacian eigenvalues of graphs, and at the end
they stated three open problems and two conjectures. The following
nice conjecture is one of them.

\begin{conjecture}\label{conj}  Let $G$ be a connected graph on $n$
vertices and $m$ edges. Then
$$q_n(G)\ge \frac{2m}{n-1}-n+2.$$\end{conjecture}

In this paper we prove a slightly stronger result:

\begin{theorem}\label{thm} Let $n$ be an integer with $n\ge 6$. If $G$ is a graph on $n$
vertices and $m$ edges, then
$$q_n(G)\ge \frac{2m}{n-2}-n+1.$$
\end{theorem}

\begin{remark} It is easy to see that $q_2(K_2)=0$. For $n\ge 3$
and $m=\frac 12n(n-1)-1$, we have $q_n(G)=\frac 12 (3n-6-\sqrt{(n-2)(n+6)})\triangleq \tau_n$.
One may find that $\tau_n<2m/(n-1)-n+2<2m/(n-2)-n+1$ for $n=3,4$
and $2m/(n-1)-n+2<\tau_n<2m/(n-2)-n+1$ for $n=5$. Thus, Conjecture
\ref{conj}  does not hold for $2\le n\le 4$, and Theorem \ref{thm}
does not hold for $2\le n\le 5$. By calculation, we find that
Conjecture \ref{conj}   holds for $n=5$.
\end{remark}

\section{Preliminary lemmas}

Given a Hermitian matrix $A$ of order $n$, we index its
eigenvalues as $\lambda_1(A)\ge \lambda_2(A)\ge \cdots \ge
\lambda_n(A)$.

\begin{lemma} (\cite{HornJohnson}) \label{lem4} 
Let $A$ and $B$ be Hermitian matrices of order $n$. Then
$$\lambda_{n}(A + B) \le \lambda_1(A) +\lambda_n(B).$$
\end{lemma}

The following Lemma \ref{lem1} is due to Merris \cite{Merris} who
stated for Laplacian eigenvalues, but his proof is also true for
signless Laplacian eigenvalues. Feng and Yu \cite{FY} gave also a
proof of Lemma \ref{lem1}.

\begin{lemma} (Merris \cite{Merris})\label{lem1}
 Let $G$ be a graph. Then
 $$q_1(G)\le\max_{u\in V(G)} \left( d(u)+\frac{1}{d(u)}\sum_{vu\in E(G)} d(v) \right).$$
\end{lemma}

\begin{lemma}{ (\cite{OLAK}, \cite{LOA})}\label{lem2} 
 Let $G\neq K_n$ be a graph of order $n$, with $k$ vertices of degree $n-1$. Then
$$q_n(G)\ge \frac{1}{2}( n+2k-2-\sqrt{(n+2k-2)^2-8k(k-1)} ).$$
\end{lemma}

\begin{lemma}{(\cite{OLAK}, \cite{LOA})}\label{lem3} 
Let $G$ be a graph of order $n$, with $k = 1$ or $2$ vertices of degree $n-1$ and at least one vertex
with degree $n-2$. If $n\ge 7$, then
$$q_n(G)\ge \frac{2k}{n-2}.$$
\end{lemma}

\section{Proof of Theorem \ref{thm}}

One may assume that $m\ge \frac 12(n-1)(n-2)+1$.  For $n=6,7$, by
Matlab, it is easy to see that Theorem 1.2 holds. Next, we assume
that $n\ge 8$. Let $m=\frac 12(n-1)(n-2)+r$. Then $r$ is a
positive integer with $r\le n-1$. Suppose that $G$ has exactly $k$
vertices of degree $n-1$. Then $G^c$ has exactly $k$ isolated
vertices, where $G^c$ is the complement of the graph $G$. It is
easy to see that Theorem \ref{thm} is equivalent to $q_n(G)\ge
{2r}/{n-2}$. Suppose that
\begin{equation}\label{eqn0}q_n(G)< \frac{2r}{n-2}.\end{equation}
We shall derive a contradiction.

By Lemma \ref{lem4} and \eqref{eqn0}, we have
\begin{equation}\label{eqn01}q_1(G^c)\ge q_n(K_n)-q_n(G)>n-2-\frac{2r}{n-2}\ge n-2-\frac r3.\end{equation}

$\blacksquare$ We shall employ Lemma \ref{lem1} to prove $r=1,2$.

By \eqref{eqn01}, $r\le n-1$ and $n\ge 8$, we have $q_1(G^c)>
n-2-(n-1)/3>0$. It follows that $E(G^c)\not=\emptyset$. By Lemma
\ref{lem1}, there exists $u\in V(G^c)$ such that
\begin{equation}\label{eqn1} q_1(G^c)\le d(u)+\frac 1{d(u)}\sum_{vu\in E(G^c)}
d(v).\end{equation} Let $v'\in V(G^c)$ with $d(v')=\max \{ d(v) :
uv\in E(G^c)\}$. Then, by \eqref{eqn1}, we have
\begin{equation}\label{eqn2} q_1(G^c)\le d(u)+d(v')\le |E(G^c)|+1= \frac 12n(n-1)-m+1=n-r.\end{equation}
It follows from  \eqref{eqn01} and \eqref{eqn2} that
$n-2-r/3<n-r$. Thus $r<3$. By \eqref{eqn01} and $r<3$, we have
$q_1(G^c)>n-3$. It follows from \eqref{eqn2} that
$d(u)+d(v')>n-3$. Hence $d(u)+d(v')\ge n-2$.

$\blacksquare$ We shall employ Lemmas \ref{lem2} and \ref{lem3} to
prove  $0\le k<r$.

If $k\ge r+1$, then, by Lemma \ref{lem2} and $n\ge 8$, we have
\begin{eqnarray*}q_n(G)&\ge&\frac{4k(k-1)}{n+2k-2+\sqrt{(n+2k-2)^2-8k(k-1)}}\\
&>&\frac{2k(k-1)}{n+2k-2}\ge \frac {4(k-1)}{n+2}\ge
\frac{4r}{n+2}>\frac{2r}{n-2},\end{eqnarray*} a contradiction with
\eqref{eqn0}. If $k=r$, then
$$\sum_{v\in V(G)} d(v) =2m = (n-1)(n-2)+2k>k(n-1)+(n-k)(n-3).$$
It follows that there exists $v\in V(G)$ such that $d(v)=n-2$.
Noting that $r=1,2$, by Lemma \ref{lem3} and $n\ge 8$, we have
$$q_n(G)\ge \frac{2k}{n-2}=\frac{2r}{n-2},$$  a contradiction with
\eqref{eqn0}. Hence $0\le k<r$.

$\blacksquare$ Exclude $r=2$.

If $r=2$, then, by $n-2\le d(u)+d(v')\le |E(G^c)|+1= n-2$, we have
$d(u)+d(v')-1=|E(G^c)|=n-3$. This implies that $G^c$ has at least
two isolated vertices, a contradiction with $k<r=2$.

$\blacksquare$ Exclude $r=1$.

Suppose that $r=1$. By \eqref{eqn0} and $n\ge 8$, we have
$q_n(G)<2r/(n-2)<1$. We shall prove that $G^c=K_{1,n-3}\cup K_2$.
Then by the eigenvalue equation of $Q(G)$ we derive  a
contradiction.

By $0\le k<r$ we have $k=0$. Thus $|E(G^c)|=n-r-1=n-2$ and $G^c$
has no isolated vertex ($k=0$). These imply that $G^c$ has at
least two nontrivial components. Let $G_u$ be the component of
$G^c$ with $u\in G_u$. So $|E(G_u)|\le |E(G^c)|-1=n-3$. It follows
from $n-2\le d(u)+d(v')\le |E(G_u)|+1\le n-2$  that
$d(u)+d(v')=n-2=|E(G_u)|+1$. This implies that, if $uv\in V(G)
(v\not= v')$, then $d(v)=1$.  If $2\le d(u)\le n-4$, then, by
\eqref{eqn1}, we have
\begin{eqnarray*} q_1(G^c)&\le& d(u)+\frac 1{d(u)}\sum_{vu\in E(G^c)}
d(v)\\
&=&d(u)+\frac 1{d(u)} (d(v')+d(u)-1)=d(u)+\frac{n-3}{d(u)}\\
&\le & n-4+\frac{n-3}{n-4}<n-2-\frac{2r}{n-2},\end{eqnarray*} a
contradiction with \eqref{eqn01}. So either $d(u)=1$, $d(v')=n-3$,
or $d(u)=n-3$, $d(v')=1$. By $|E(G_u)|=n-3$, we have $G_u=K_{1,\,
n-3}$ and then $G^c=K_{1,\, n-3}\cup K_2$.

Denote by $v_1$ the center of the star $K_{1,\, n-3}$, by $v_2,
v_3$ the end vertices of $K_2$, and by $v_4,\ldots,v_n$ the
pendent vertices of the star $K_{1,\, n-3}$ respectively. Then $G$
is a connected graph with vertex set $V(G)=\{\, v_1, v_2, \ldots,
v_n\}$.
 Let $X=(x_1, x_2, \ldots, x_n)^T$ be an eigenvector of $Q(G)$ corresponding to
 $q_n(G)$, i.e. $Q(G)X=q_n(G)X$.  By the symmetry of $G$ and $q_n(G)<1$,
  we have $x_2=x_3$, $x_4=\cdots =x_n$. From the eigenvalue equation, we have
\begin{eqnarray*}
&&(q_n(G)-2)x_1-2x_2=0,\\
&&2x_2+(-q_n(G)+2n-6)x_4=0,\\
&&x_1+(-q_n(G)+n-2)x_2+(n-3)x_4=0.
\end{eqnarray*}
Since at least one of $x_1, x_2, x_4$ is nonzero, it follows from
the above system that $x_1, x_2, x_4$ are all nonzero. Noting that
$0\le q_n(G)<1$ and $n\ge 8$, we have
$$|x_1|<|(q_n(G)-2)x_1|=2|x_2|,$$
$$(n-3)|x_4|<(-q_n(G)+2n-6)|x_4|=2|x_2|,$$
$$(n-3)|x_2|<(-q_n(G)+n-2)|x_2|=|x_1+(n-3)x_4|<4|x_2|,$$ a
contradiction with $n\ge 8$.

  This completes the
proof of Theorem \ref{thm}.

\vskip 2mm \noindent{\bf Acknowledgments }

\vskip 3mm We are grateful to  the anonymous referees for valuable
suggestions  which result in an improvement of the original
manuscript.

\small {

}

\end{document}